\documentclass[11pt]{article}
\usepackage{}

\usepackage[total={6in, 8.5in}]{geometry}
\usepackage{amsfonts}
\usepackage{amsthm}
\usepackage{amssymb}
\usepackage{amsmath}
\usepackage{enumerate}
\usepackage[
pdfauthor={},
pdftitle={},
pdfstartview=XYZ,
bookmarks=true,
colorlinks=true,
linkcolor=blue,
urlcolor=blue,
citecolor=blue,
bookmarks=false,
linktocpage=true,
hyperindex=true
]{hyperref}

\makeatletter
\newcommand{\setword}[2]{%
	\phantomsection
	#1\def\@currentlabel{\unexpanded{#1}}\label{#2}%
}
\makeatother

\makeatletter
\renewcommand*\env@matrix[1][*\c@MaxMatrixCols c]{%
	\hskip -\arraycolsep
	\let\@ifnextchar\new@ifnextchar
	\array{#1}}
\makeatother

\usepackage{graphicx}
\usepackage[natural]{xcolor}
\usepackage{tikz}
\usepackage{tikz-3dplot}

\usepackage{pgf,tikz,pgfplots}

\usepackage{mathrsfs}

\usetikzlibrary{arrows}

\usepackage{xcolor}
\long\def\ignore#1{}

\newtheorem{THM}{Theorem}

\newtheorem{LEM}[THM]{Lemma}
\newtheorem{CON}{Conjecture}

\newtheorem{COR}[THM]{Corollary}

\newtheorem{CLA}{Claim}
\newcommand{\pf}{\textbf{Proof}.\quad}

\newtheorem*{LEM3}{\textbf{Lemma 5}}

\usepackage{thmtools}
\usepackage{thm-restate}

\newcommand{\CC}{\mathcal{C}}

\newcommand{\pbar}{\overline{\varphi}}
\DeclareMathOperator{\dist}{dist} 

\begin{document}
\title{An improvement to the vertex-splitting conjecture}

\author{%
 Yan Cao
 \qquad Guantao Chen\thanks{This author was supported in part by NSF grant DMS-1855716.}\\
  Department of Mathematics and Statistics, \\
  Georgia State University, Atlanta, GA 30302, USA\\
   \texttt{ycao17@gsu.edu}
  \qquad
   \texttt{gchen@gsu.edu}%
 \and 
 Songling Shan \\
 Department of Mathematics, \\
 Illinois State Univeristy, Normal, IL 61790, USA \\
 \texttt{sshan12@ilstu.edu}
} 

\date{March 8, 2021}
\maketitle

 \begin{abstract}
 	
For a simple graph $G$,  denote by $n$, $\Delta(G)$,
and $\chi'(G)$ its order, maximum degree, and chromatic index, 
respectively. A connected class 2 graph $G$ is \emph{edge-chromatic critical} 
if $\chi'(G-e)<\Delta(G)+1$ for every edge $e$ of $G$. 
Define $G$ to be   \emph{overfull}  if $|E(G)|>\Delta(G) \lfloor n/2 \rfloor$.
Clearly, overfull graphs are class 2 and any graph obtained from a regular graph of even order by splitting a vertex is overfull. 
Let $G$ be an $n$-vertex connected regular class 1 graph with $\Delta(G) >n/3$.
Hilton and Zhao in  1997 conjectured that if $G^*$ is obtained from $G$
by splitting one vertex of $G$ into two vertices, then $G^*$
is edge-chromatic critical, and  they verified the conjecture for graphs $G$ 
with $\Delta(G)\ge \frac{n}{2}(\sqrt{7}-1)\approx 0.82n$. 
The graph $G^*$ is easily verified to be overfull, and so the hardness of the 
conjecture lies in showing that the deletion of every of its edge decreases the chromatic index. Except in 2002, Song showed that the conjecture is true for 
a special class of graphs $G$ with $\Delta(G)\ge \frac{n}{2}$,
no other progress on this conjecture had been made. 
In this paper, we confirm the conjecture for graphs $G$ with  $\Delta(G) \ge  0.75n$.


 \end{abstract}

\section{Introduction}

We consider only simple graphs. 
  Let $G$ be a graph.  Denote by $V(G)$ and $E(G)$ the vertex set and edge set of $G$, respectively. 
  An  \emph{edge $k$-coloring}  of $G$ is a mapping $\varphi$ from $E(G)$ to the set of integers
  $[1,k]:=\{1,\ldots, k\}$, called {\it colors\/}, such that  no two adjacent edges receive the same color with respect to $\varphi$.  
  The {\it chromatic index\/} of $G$, denoted $\chi'(G)$, is defined to be the smallest integer $k$ so that $G$ has an edge $k$-coloring.  
  We denote by $\CC^k(G)$ the set of all edge $k$-colorings of $G$.  In 1960's, Vizing~\cite{Vizing-2-classes} showed that every simple graph  $G$ has  chromatic index either $\Delta(G)$ or $\Delta(G)+1$.
  If $\chi'(G)=\Delta(G)$, then $G$ is said to be of {\it class 1\/}; otherwise, it is said to be
  of {\it class 2\/}.  
  Holyer~\cite{Holyer} showed that it is NP-complete to determine whether an arbitrary graph is of class 1.
  However, if a graph $G$ has too many edges, i.e., $|E(G)|>\Delta(G) \lfloor |V(G)|/2\rfloor$,  then $G$
  is class 2. Such graphs are called \emph{overfull}. 
 Easily implied by its definition, overfull graphs are of odd order.

 We call $G$
  \emph{edge-chromatic critical} or \emph{$\Delta$-critical} if $\chi'(G)=\Delta(G)$+1 and $\chi'(H)<\Delta(G)+1$ for every proper subgraph $H$ of $G$. For example, odd cycles and the graph obtained from the Petersen graph by deleting one vertex are edge-chromatic critical.  We study sufficient conditions for a class 2 graph to be edge-chromatic critical. 
    A \emph{vertex-splitting} in $G$ at a vertex $v$ is obtained 
    by replacing $v$ with two new adjacent vertices $v_1$
    and $v_2$ and partition the neighborhood $N_G(v)$
    into two nonempty subsets that serve as the neighborhoods of $v_1$ and $v_2$ in $G'$, respectively. 
     We say $G'$ is obtained from $G$ by a vertex-splitting. 
  A vertex-splitting was formulated in terms of the ``M\H{o}bius-type gluing technique''
  in~\cite{MR4183149} and~\cite{MR3417245}. 
  Hilton and Zhao~\cite{MR1460574} in 1997 propsed the following conjecture. 
  \begin{CON}[Vertex-splitting conjecture]
   Let $G$ be an $n$-vertex class 1 $\Delta$-regular graph with $\Delta>\frac{n}{3}$. If $G^*$
  is obtained from $G$ by a \emph{vertex-splitting}, 
  then $G^*$ is $\Delta$-critical. 
  \end{CON} 

Since the graph $G^*$ above is overfull and so is class 2, the difficulty of the vertex-splitting conjecture lies in checking every edge of $G^*$
is critical, i.e., whose deletion decreases the chromatic index of $G^*$.   
 Hilton and Zhao~\cite{MR1460574} in the same paper verified the conjecture for graphs $G$ 
  with $\Delta(G)\ge \frac{n}{2}(\sqrt{7}-1)\approx 0.82n$. 
 Song~\cite{MR1874750} in 2002 showed that the conjecture is true for 
  a special class of graphs $G$ with $\Delta(G)\ge \frac{n}{2}$.
  Except this result, to our best knowledge, we are not aware of 
any other progress on the conjecture. 
  In this paper, we verify the conjecture for graphs $G$ with  $\Delta(G) \ge  0.75n$ as below.

  \begin{THM}\label{thm:vertex spliting}
  	Let $n$ and $ \Delta$ be positive integers such that $\Delta\ge \frac{3(n-1)}{4}$. 
  	If $G$ is obtained from an $(n-1)$-vertex $\Delta$-regular class 1 graph by a vertex-splitting, then $G$ 
  	is $\Delta$-critical. 
  \end{THM}

The reminder of this paper is organized as follows. In Section 2, 
we introduce some definitions and preliminary results.  
In Section 3, we prove Theorem~\ref{thm:vertex spliting}. In the last 
Section, we prove one newly developed adjacency lemma.

\section{Definitions and Preliminary Results}\label{lemma}

Let $G$ be a graph. 
For $e\in E(G)$, $G-e$
denotes the graph obtained from $G$ by deleting the edge $e$. 
The symbol $\Delta$  is reserved for $\Delta(G)$, the maximum degree of $G$
throughout  this paper.  A \emph{$k$-vertex}  in $G$  is a vertex of degree  $k$
in $G$, and a \emph{$k$-neighbor}  of a vertex $v$ is a neighbor of $v$ that is a $k$-vertex in $G$.
For $u,v\in V(G)$, we use $\dist_G(u,v)$ to denote the distance between $u$ and $v$, which is the length of a shortest path connecting $u$
and $v$ in $G$. For $S\subseteq V(G)$, define $\dist_G(u,S)=\min_{v\in S} \dist_G(u,v)$. 

An edge $e\in E(G)$ is a \emph{critical edge} of $G$ if $\chi'(G-e)<\chi'(G)$. 
It is not hard to see that for a connected class 2 graph, if 
every of its edge   is critical, then $G$ is $\Delta$-critical.  
Edge-chromatic critical graphs provide more information about the  structure around a vertex than general class 2 graphs. For 
example,  Vizing's Adjacency Lemma (VAL) from 1965~\cite{Vizing-2-classes} is a useful tool that reveals certain structure at a vertex
by assuming the criticality of an edge.

\begin{LEM}[Vizing's Adjacency Lemma (VAL),\cite{Vizing-2-classes}]Let $G$ be a class 2 graph with maximum degree $\Delta$. If $e=xy$ is a critical edge of $G$, then $x$ has at least $\Delta-d_G(y)+1$ $\Delta$-neighbors in $V(G)\setminus \{y\}$.
	\label{thm:val}
\end{LEM}
Let $G$ be a graph and 
$\varphi\in \CC^k(G-e)$ for some edge $e\in E(G)$ and some integer $k\ge 0$. 
For any $v\in V(G)$, the set of colors \emph{present} at $v$ is 
$\varphi(v)=\{\varphi(f)\,:\, \text{$f$ is incident to $v$}\}$, and the set of colors \emph{missing} at $v$ is $\pbar(v)=[1,k]\setminus\varphi(v)$.  
For a vertex set $X\subseteq V(G)$,  define 
$
\pbar(X)=\bigcup _{v\in X} \pbar(v).
$
We call $X$  \emph{elementary} with respect to $\varphi$  or simply \emph{$\varphi$-elementary} if $\pbar(u)\cap \pbar(v)=\emptyset$
for every two distinct vertices $u,v\in X$.   Sometimes, we just say that $X$ 
is elementary if the  edge coloring is understood.

 For two distinct colors $\alpha,\beta \in [1,k]$, let  $H$ be the spanning subgraph of $G$
 with its edges  from $E(G)$ that are colored by $\alpha$
 or $\beta$ with respect to $\varphi$. Each component of $H$ is either 
 an even cycle or a path, which is called an \emph{$(\alpha,\beta)$-chain} of $G$
 with respect to $\varphi$.  If we interchange the colors $\alpha$ and $\beta$
 on an $(\alpha,\beta)$-chain $C$ of $G$, we get a new edge $k$-coloring  of $G$, 
 and we write $$\varphi'=\varphi/C.$$
 This operation is called a \emph{Kempe change}.

 Let  $x,y\in V(G)$, and  $\alpha, \beta\in [1,k]$ be two distinct colors.   If $x$ and $y$
 are contained in a same  $(\alpha,\beta)$-chain of $G$ with respect to $\varphi$, we say $x$ 
 and $y$ are \emph{$(\alpha,\beta)$-linked} with respect to $\varphi$.
 Otherwise, $x$ and $y$ are \emph{$(\alpha,\beta)$-unlinked} with respect to $\varphi$. For a vertex-edge sequence $S$, we use $V(S)$
 to denote the set of all vertices contained in the sequence. 
 
\subsection{Multifan and Kierstead Path}
The fan argument was introduced by Vizing~\cite{Vizing64,vizing-2factor} in his classical results on the upper bounds of chromatic indices for simple graphs and multigraphs.  We will use multifan, a generalized version of Vizing fan, given by Stiebitz et al.~\cite{StiebSTF-Book}, in our proof.

Let  $G$ be a graph, $e=rs_1\in E(G)$ and $\varphi\in \CC^k(G-e)$ for some integer $k\ge 0$.
A \emph{multifan} centered at $r$ with respect to $e$ and $\varphi$
is a sequence $F_\varphi(r,s_1:s_p):=(r, rs_1, s_1, rs_2, s_2, \ldots, rs_p, s_p)$ with $p\geq 1$ consisting of  distinct vertices $r, s_1,s_2, \ldots , s_p$ and distinct edges $rs_1, rs_2,\ldots, rs_p$ satisfying the following condition:
\begin{enumerate}[(F1)]
	\item For every edge $rs_i$ with $i\in [2,p]$,  there exists $j\in [1,i-1]$ such that 
	$\varphi(rs_i)\in \pbar(s_j)$. 
\end{enumerate}
We will simply denote a multifan  $F_\varphi(r,s_1: s_{p})$ by $F$ if 
$\varphi$ and the vertices and edges in $F_\varphi(r,s_1: s_{p})$  are clear. 
The following result regarding a multifan can be found in \cite[Theorem~2.1]{StiebSTF-Book}.

\begin{LEM}
	\label{thm:vizing-fan1}
	Let $G$ be a class 2 graph and $F_\varphi(r,s_1:s_p)$  be a multifan with respect to a critical edge $e=rs_1$ and a coloring $\varphi\in \CC^\Delta(G-e)$. Then  the following statements  hold. 
	 \begin{enumerate}[(a)]
	 	\item $V(F)$ is $\varphi$-elementary. \label{thm:vizing-fan1a}
	 	\item Let $\alpha\in \pbar(r)$. Then for every $i\in [1,p]$  and $\beta\in \pbar(s_i)$,  $r$ 
	 	and $s_i$ are $(\alpha,\beta)$-linked with respect to $\varphi$. \label{thm:vizing-fan1b}
	 \end{enumerate}
\end{LEM}

Let $G$ be a graph, $e=v_0v_1\in E(G)$, and  $\varphi\in \CC^k(G-e)$ for some integer $k\ge 0$.
A \emph{Kierstead path}  with respect to $e$ and $\varphi$
is a sequence $K=(v_0, v_0v_1, v_1, v_1v_2, v_2, \ldots, v_{p-1}, v_{p-1}v_p,  v_p)$ with $p\geq 1$ consisting of  distinct vertices $v_0,v_1, \ldots , v_p$ and distinct edges $v_0v_1, v_1v_2,\ldots, v_{p-1}v_p$ satisfying the following condition:
\begin{enumerate}[(K1)]
	\item For every edge $v_{i}v_{i+1}$ with $i\in [1,p-1]$,  there exists $j\in [0,i-1]$ such that 
	$\varphi(v_{i}v_{i+1})\in \pbar(v_j)$. 
\end{enumerate}

Clearly a Kierstead path with at most 3 vertices is a multifan. We consider Kierstead paths with $4$ vertices. The result below was proved in Theorem 3.3 from~\cite{StiebSTF-Book}. 

\begin{LEM}[]\label{Lemma:kierstead path1}
	Let $G$ be a class 2 graph,
	 $e=v_0v_1\in E(G)$ be a critical edge, and $\varphi\in \CC^\Delta(G-e)$. If $K=(v_0, v_0v_1, v_1, v_1v_2,  v_2, v_2v_3, v_3)$ is a Kierstead path with respect to $e$
		and $\varphi$, then the following statements hold.
	\begin{enumerate}[(a)]
 		\item If $\min\{d_G(v_1), d_G(v_2)\}<\Delta$, then $V(K)$ is $\varphi$-elementary.
 		\item $|\pbar(v_3)\cap (\pbar(v_0)\cup \pbar(v_1))|\le 1$. 
 	\end{enumerate}

\end{LEM}

	\section{Proof of Theorem~\ref{thm:vertex spliting}}

\begin{figure}[h]
	\centering 
	\includegraphics[scale=0.6]{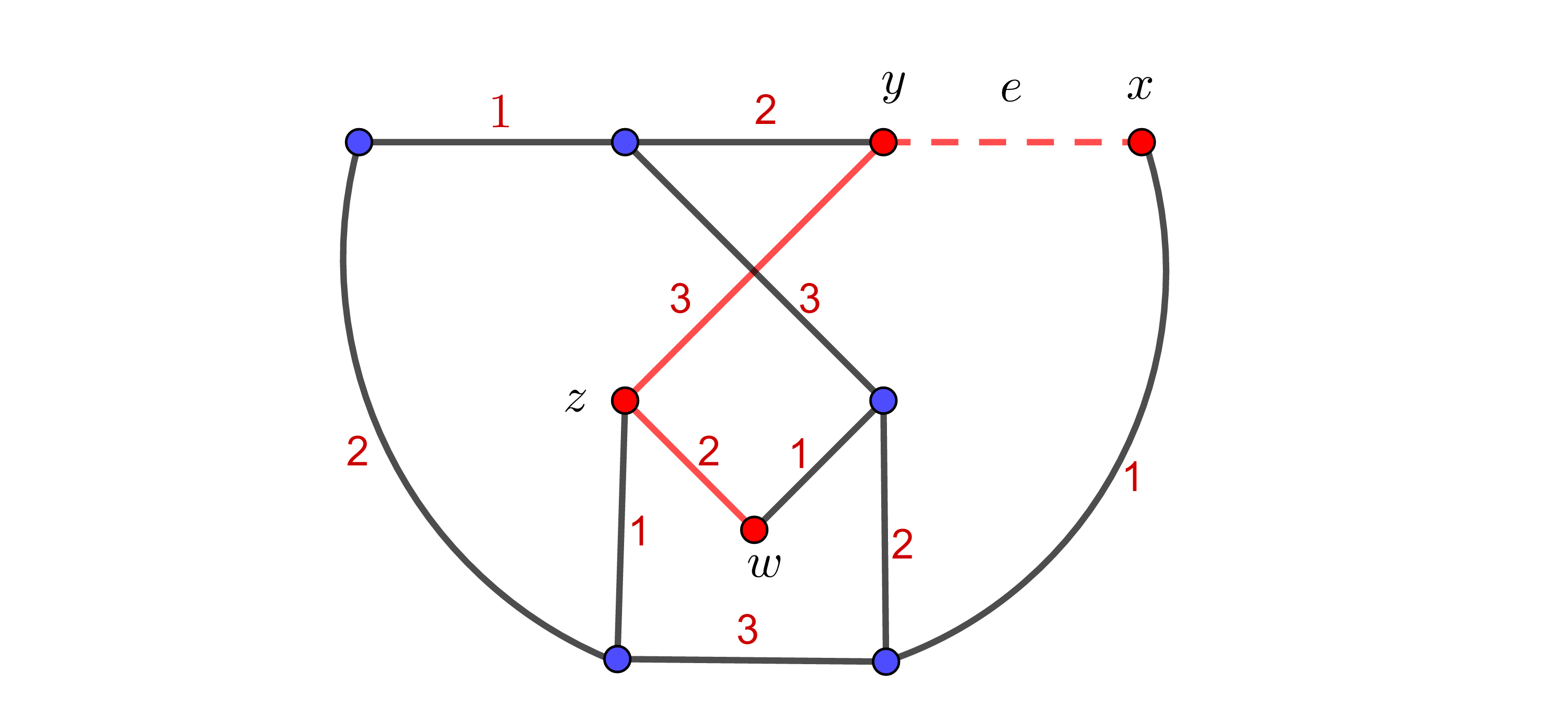}\\
	\caption{A  Kierstead path with non-elementary vertex set in a 3-coloring   of $P^*-xy$}
	\label{figP}
\end{figure}

The proof of Theorems~\ref{thm:vertex spliting} is mainly an application 
of a new adjacency lemma--Lemma~\ref{lemma:class2-with-fullDpair2} below. We define 
a {\it short-kite}  to be a 6-vertex graph consisting of a 4-cycle $abuca$
and two additional edges $ux$ and $uy$.  The truth of the vertex-splitting conjecture would be evident when $\Delta \ge \frac{n}{2}$ if the vertex set of every 
Kierstead path on four vertices is elementary. 
Unfortunately, the statement is not true and a counterexample has been 
found, see Figure~\ref{figP}, where $K=(x,xy,y,yz,z,zw,w)$ and $V(K)$
is not elementary with respect to the given coloring, and $P^*$
is obtained from the Petersen graph by deleting one vertex.

The new adjacency lemma below 
is an attempt to reveal some elementary properties of a Kierstead path on four vertices  by incorporating some additional structure to the path. 

\begin{LEM}\label{lemma:class2-with-fullDpair2}
	Let  $G$ be  a class 2 graph, 
	 $H\subseteq G$ 
	be a short-kite with $V(H)=\{a,b,c,u,x,y\}$, and let $\varphi\in \CC^\Delta(G-ab)$. 
	Suppose $$K=(a,ab,b,bu,u,ux,x) \quad \text{and} \quad K^*=(b,ab,a,ac,c,cu,u,uy,y)$$
	are two Kierstead path with respect to $ab$ and $\varphi$.  
	If $\pbar(x)\cup \pbar(y)\subseteq \pbar(a)\cup \pbar(b)$,  then $\max\{d_G(x),d_G(y)\}=\Delta $. 
\end{LEM}

The proof of Lemma~\ref{lemma:class2-with-fullDpair2} will be  given in the last section. 
Since all vertices not missing  a given color $\alpha$
are saturated by the matching that consists of all edges colored by $\alpha$ in $G$, we have the Parity Lemma below, which has appeared in many papers, for example, see~\cite[Lemma 2.1]{MR2028248}.

\begin{LEM}[Parity Lemma]
	Let $G$ be an $n$-vertex graph and $\varphi\in \CC^\Delta(G)$. 
	Then for any color $\alpha\in [1,\Delta]$, 
	$|\{v\in V(G): \alpha\in \pbar(v)\}| \equiv n \pmod{2}$. 
\end{LEM}

Let $G$ be a graph and $u,v\in V(G)$ be adjacent. 
We call $(u,v)$ a \emph{full-deficiency pair} of $G$
if $d(u)+d(v)=\Delta(G)+2$.  If $G$ is $\Delta$-critical, then a full-deficiency pair $(u,v)$ of $G$ is called a \emph{saturating pair} of $G-uv$ in~\cite{MR4183149}.  
%
%

\begin{LEM}\label{lemma:class2-with-fullDpair}
	If  $G$ is an $n$-vertex class 2 graph with a full-deficiency pair $(a,b)$ such that $ab$ is a critical edge of $G$, 
	then $G$ satisfies the following properties. 
	\begin{enumerate}[$(i)$]
		
		\item For every $x\in (N_G(a)\cup N_G(b))\setminus\{a,b\}$, $d_G(x)=\Delta$;
		
		\item For every $x\in V(G)\setminus\{a,b\}$, if $\dist_G(x, \{a,b\})=2$, then 
		$d_G(x)\ge \Delta-1$. 
		Furthermore, 
		if $d_G(a)<\Delta$ and $d_G(b)<\Delta$, then  $d_G(x)=\Delta$;
		\item  For  
		every   $x\in V(G)\setminus\{a,b\}$,  if $d_G(x)\ge n-|N_G(b)\cup N_G(a)|$, 
		then $d_G(x)\ge \Delta-1$. 
		Furthermore, 
			if $d_G(a)<\Delta$ and $d_G(b)<\Delta$, then  $d_G(x)=\Delta$;
		\item Suppose that $n$ is odd. If there exists $x\in V(G)\setminus \{a,b\}$ such that $d_G(x)<\Delta$, then there exists 
		$y\in V(G)\setminus \{a,b,x\}$ such that $d_G(y)<\Delta$.
		%
	\end{enumerate} 
\end{LEM}
 \pf We let
 $\varphi\in \CC^\Delta(G-ab)$ and 
 $
 F=(b, ba,a) 
 $
 be the multifan  with respect to $ab$ and $\varphi$. 
 By Lemma~\ref{thm:vizing-fan1}~\eqref{thm:vizing-fan1a}, 
 \begin{equation}\label{pbarFa}
 |\pbar(V(F))|= 2\Delta-(d_{G}(a)+d_{G}(b)-2)= 2\Delta+2-(\Delta+2)=\Delta.
 \end{equation}
 
 By Lemma~\ref{thm:vizing-fan1}, for every $\varphi'\in \CC^\Delta(G-ab)$,  $\{a,b\}$ is $\varphi'$-elementary and for every 
 $i\in \pbar'(a)$ and $j\in \pbar'(b)$, $a$ and $b$ are $(i,j)$-linked with respect to $\varphi'$. We will use this fact very often.

Since $\pbar(a)\cap \pbar(b)=\emptyset$ and $\pbar(a)\cup \pbar(b)=[1,\Delta]$, it follows that $\varphi(a)=\pbar(b)$. 
Thus for any $x\in N_G(a)\setminus \{b\}$, $(a, ab,b,ax,x)$ 
is a multifan with respect to $ab$ and $\varphi$ and so $\{a,b,x\}$ is $\varphi$-elementary by Lemma~\ref{thm:vizing-fan1} (a).   It follows from~\eqref{pbarFa}
that $d_G(x)=\Delta$. Symmetrically, for each  $x\in N_G(b)\setminus \{a\}$, $d_G(x)=\Delta$. 
This proves (i).

For (ii), let $x\in V(G)\setminus\{a,b\}$ such that $\dist_G(x, \{a,b\})=2$. 
We assume that $\dist_G(x, b)=2$ and let 
$u\in ( N_G(b))\setminus\{a\})\cap N_G(x)$.
Then by~\eqref{pbarFa}, $K=(a,ab, b, bu, u, ux, x)$ is a Kierstead path with respect to $ab$
and $\varphi$. By~\eqref{pbarFa} and Lemma~\ref{Lemma:kierstead path1} (b),  it follows that 
$d_G(x)\ge \Delta-1$. If $d_G(a)<\Delta$ and $d_G(b)<\Delta$, then
$V(K)$ is $\varphi$-elementary  by  Lemma~\ref{Lemma:kierstead path1} (a).
Since $\pbar(a)\cup \pbar(b)=[1,\Delta]$ by~\eqref{pbarFa},
it follows that  $d_G(x)=\Delta$.

For (iii), let $x\in V(G)\setminus\{a,b\}$ such that $d_G(x)\ge n-|N_G(b)\cup N_G(a)|$. By (i), we may assume that 
$x\not\in (N_G(a)\cup N_G(b))\setminus\{a,b\}$. Thus   $d_G(x)\ge n-|N_G(b)\cup N_G(a)|$ implies that 
there exists $u\in ( (N_G(a)\cup N_G(b)))\cap N_G(x)$. 
Therefore, $\dist_G(x,\{a,b\})=2$. Now Statement (ii) yields the conclusion.
For any color $\alpha\in \pbar(x)$,  either $\alpha\in \pbar(a)$ or $\alpha\in \pbar(b)$, since $\{a,b\}$ is $\varphi$-elementary and $\pbar(a)\cup \pbar(b)=[1,\Delta]$. Since $n$ is odd, 
by the Parity Lemma, there exists  
$y\in V(G)\setminus \{a,b,x\}$ such that $\alpha\in \pbar(y)$,
and so $d_G(y)<\Delta$, proving (iv).  
%
\qed 

\begin{COR}\label{cor:no2D-1}
	Let  $G$ be  an $n$-vertex class 2 graph with a full-deficiency pair $(a,b)$ such that $ab$ is a critical edge of $G$.  
	If $\Delta\ge \frac{3(n-1)}{4}$, then there exists at most one vertex $x\in V(G)\setminus \{a,b\}$ such that $d_G(x)=\Delta-1$.
\end{COR}	

\pf Assume to the contrary that there exist distinct $x,y\in V(G)\setminus \{a,b\}$ 
such that $d_G(x)=d_G(y)=\Delta-1$. By Lemma~\ref{lemma:class2-with-fullDpair} (i), $x,y\not\in (N_G(a)\cup N_G(b))\setminus\{a,b\}$. 
By Lemma~\ref{lemma:class2-with-fullDpair} (iii), we may assume  that $d_G(b)=\Delta$. 
Thus $d_G(a)=2$ as $d_G(a)+d_G(b)=\Delta+2$. 
Let $c$ be the other neighbor of $a$ in $G$. Since $(a,c)$ is a full-deficiency pair of $G$ as well,  we may assume $x,y\not\in N_G(c)$.

Since $d_{G}(b)=d_{G}(c)=\Delta$ and  $d_{G}(x)=d_{G}(y)=\Delta-1$, we get $|N_{G}(b)\cap N_{G}(c)|\ge \frac{n}{2}-1$ and $|N_{G}(x)\cap N_{G}(y)|\ge \frac{n}{2}-2$. 
Since $b,c,x,y\not\in N_{G}(b)\cap N_{G}(c)$ and $b,c,x,y\not\in N_{G}(x)\cap N_{G}(y)$, we get 
$|N_{G}(b)\cap N_{G}(c)\cap N_{G}(x)\cap N_G(y)|\ge 1$.
Let 
$u\in N_{G}(b)\cap N_{G}(c)\cap N_{G}(x)\cap N_{G}(y)$,  $H$ be 
the short-kite with $V(H)=\{a,b,c,u,x,y\}$, and $\varphi\in \CC^\Delta(G-ab)$. 
As $\{a,b\}$
is $\varphi$-elementary$, |\pbar(a)\cup \pbar(b)|=2\Delta+2-(d_{G}(a)+d_{G}(b))= \Delta$ and so $\pbar(a)\cup \pbar(b)=[1,\Delta]$. Thus   
$K=(a,ab,b,bu,u,ux,x)$ and $ K^*=(b,ab,a,ac,c,cu,u,uy)$
are two Kierstead paths with respect to $ab$ and $\varphi$, and $\pbar(x)\cup \pbar(y)\subseteq \pbar(a)\cup \pbar(b)$.   
However,  $d_{G}(x)=d_{G}(y)=\Delta-1$,  contradicting  Lemma~\ref{lemma:class2-with-fullDpair2}.
\qed  

\proof[\bf Proof of Theorem~\ref{thm:vertex spliting}]
Since $G$ is overfull, it is  class 2. 
We only need to show that every edge of $G$ is critical. 
Suppose to the contrary that there exists $xy\in E(G)$ such that $xy$ 
is not a critical edge of $G$.  Let 
$
G^*=G-xy.  
$
Then $\chi'(G^*)=\Delta+1$.

Since $ab$ is a critical edge of $G$,  $ab \ne xy$. 
Also, since $ab$ is a critical edge of $G$, and any $\Delta$-coloring of $G-ab$
gives a $\Delta$-coloring of $G^*-ab$, $ab$ is also a critical edge of $G^*$. 
Now $d_{G^*}(x)=d_{G^*}(y)=\Delta-1$,  reaching a contradiction to Corollary~\ref{cor:no2D-1}. 
\qed

\section{Proof of Lemma~\ref{lemma:class2-with-fullDpair2}}

 
 We start with some notation. 
 Let $G$ be a graph and 
 $\varphi\in \CC^k(G-e)$ for some edge $e\in E(G)$ and some integer $k\ge 0$. 
For all the concepts below, when we use them later on, 
if we skip $\varphi$, we mean the concept is defined with respect to the current edge coloring.

Let  $x,y\in V(G)$, and  $\alpha, \beta, \gamma\in [1,k]$ be three colors.   
Let $P$ be an 
$(\alpha,\beta)$-chain of $G$ with respect to $\varphi$ that contains both $x$ and $y$. 
If $P$ is a path, denote by $\mathit{P_{[x,y]}(\alpha,\beta, \varphi)}$  the subchain  of $P$ that has endvertices $x$
and $y$.  By \emph{swapping  colors} along  $P_{[x,y]}(\alpha,\beta,\varphi)$, we mean 
exchanging the two colors $\alpha$
and $\beta$ on the path $P_{[x,y]}(\alpha,\beta,\varphi)$.

Define  $P_x(\alpha,\beta,\varphi)$ to be an $(\alpha,\beta)$-chain or an $(\alpha,\beta)$-subchain of $G$ with respect to $\varphi$ that starts at $x$  and ends at a different vertex missing exactly one of $\alpha$ and $\beta$.    
If $x$ is an endvertex of the $(\alpha,\beta)$-chain that contains $x$, then $P_x(\alpha,\beta,\varphi)$ is unique.  Otherwise, we take one segment of the whole chain to be 
$P_x(\alpha,\beta,\varphi)$. We will specify the segment when it is used.

If  $u$  is a vertrex on  $P_x(\alpha,\beta,\varphi)$, we  write  {$\mathit {u\in P_x(\alpha,\beta, \varphi)}$}; and if $uv$  is an edge on  $P_x(\alpha,\beta,\varphi)$, we  write  {$\mathit {uv\in P_x(\alpha,\beta, \varphi)}$}.  
If $u,v\in P_x(\alpha,\beta,\varphi)$ such that $u$ lies between $x$ and $v$ on the path, 
then we say that $P_x(\alpha,\beta,\varphi)$ \emph{meets $u$ before $v$}. 
Suppose the current color of an  edge $uv$ of $G$
is $\alpha$, the notation  $\mathit{uv: \alpha\rightarrow \beta}$  means to recolor  the edge  $uv$ using the color $\beta$.  If $|\pbar(x)|=1$, we will also use $\pbar(x)$ to denote the  color that is missing at $x$.

Let $\alpha, \beta, \gamma, \tau,\eta\in[1,k]$. 
We will use a  matrix with two rows to denote a sequence of operations  taken  on $\varphi$.
Each entry in the first row represents a path or  a sequence of vertices. 
Each entry in the second row, indicates the action taken on the object above this entry. 
We require the operations to be taken to follow the ``left to right'' order as they appear in 
the matrix. 
For example,   the matrix below indicates 
three sequential operations taken on the graph based 
on the coloring from the previous step:
\[
\begin{bmatrix}
P_{[a, b]}(\alpha, \beta) &   rs & ab \\
\alpha/\beta & \gamma \rightarrow \tau & \eta
\end{bmatrix}.
\]
\begin{enumerate}[Step 1]
	\item Swap colors on the $(\alpha,\beta)$-subchain $P_{[a, b]}(\alpha, \beta,\varphi) $.
	
	\item Do  $rs: \gamma \rightarrow \tau $. 
	\item Color the edge $ab$ using color $\eta$. 
\end{enumerate}

\begin{LEM3}
	Let  $G$ be  a class 2 graph, 
	$H\subseteq G$ 
	be a short-kite with $V(H)=\{a,b,c,u,x,y\}$, and let $\varphi\in \CC^\Delta(G-ab)$. 
	Suppose $$K=(a,ab,b,bu,u,ux,x) \quad \text{and} \quad K^*=(b,ab,a,ac,c,cu,u,uy,y)$$
	are two Kierstead path with respect to $ab$ and $\varphi$.  
	If $\pbar(x)\cup \pbar(y)\subseteq \pbar(a)\cup \pbar(b)$,  then $\max\{d_G(x),d_G(y)\}=\Delta $. 
\end{LEM3}
 
 \pf 
 Assume to the contrary that  $\max\{d_G(x),d_G(y)\}\le \Delta-1$. 
 Since both $K$ and $K^*$ are Kierstead paths and 
 $\pbar(x)\cup \pbar(y)\subseteq \pbar(a)\cup \pbar(b)$, Lemma~\ref{Lemma:kierstead path1} (a)
 and (b) implies that $d_G(b)=d_G(u)=\Delta$
 and $d_G(x)=d_G(y)=\Delta-1$. 
 
 
 Let $\pbar(b)=\{1\}$. Then $\varphi(ac)=1$. 
 We may assume  $\varphi(uy)=1$.  The reasoning is below. 
 Since $a$ and $b$ are $(1,\alpha)$-linked for every $\alpha\in \pbar(y)\subseteq \pbar(a)\cup \pbar(b)$, 
 we may assume  $\pbar(y)=1$. Then a $(1,\varphi(uy))$-swap at  $y$
 gives a coloring, call it still $\varphi$, such that $\varphi(uy)=1$. 
 We consider now two cases. 
 
 \smallskip 
 {\bf \noindent Case 1: $\pbar(x)=\pbar(y)$.} 
 \smallskip 
 
 Let 
 $ \varphi(ux)=\gamma$,  \text{and}  $\pbar(x)=\pbar(y)=\eta.$	
 As $\varphi(uy)=\pbar(b)=1$, $1\not\in \{\gamma, \eta\}$. 
 As both $K$ and $K^*$ are Kierstead paths and 
 $\pbar(x)\cup \pbar(y)\subseteq \pbar(a)\cup \pbar(b)$, $\gamma,\eta \in \pbar(a)$. 
 Denote by $P_u(1,\gamma)$  the $(1,\gamma)$-subchain starting at $u$ that does not include the edge $ux$. 
 \begin{CLA}\label{cla:claim1}
 	We may assume that $P_u(1,\gamma)$    ends at $x$, some vertex  $z\in V(G)\setminus\{a,b,c,u,x,y\}$,  or passing $c$ ends at $a$. 
 \end{CLA}
 
 \pf Note that $P_a(1,\gamma)=P_b(1,\gamma)$. If $u\not\in P_a(1,\gamma)$, then the $(1,\gamma)$-chain containing $u$ is a cycle or a path with endvertices contained in $V(G)\setminus\{a,b,c,u,x,y\}$.  Thus  
 $P_u(1,\gamma)$ ends at  $x$ or some  $z\in V(G)\setminus\{a,b,c,u,x,y\}$.  Hence we assume $u\in P_a(1,\gamma)$. 
 As a consequence,  $P_u(1,\gamma)$ ends at either $b$ or $a$. 
 If $P_x(1,\gamma)$ ends at $b$, we color $ab$ by 1, uncolor $ac$, and exchange the vertex labels $b$ and $c$. 
 This gives an edge $\Delta$-coloring of $G-ab$ such that  $P_u(1,\gamma)$ ends at  $a$.
 Thus, if $u\in P_a(1,\gamma)$, we may always assume that $P_u(1,\gamma)$ ends at  $a$.	  
 \qed 
 
\begin{figure}[!htb]
 	\begin{center}
 		\begin{tikzpicture}[scale=0.8]
 		
 		{\tikzstyle{every node}=[draw ,circle,fill=white, minimum size=0.5cm,
 			inner sep=0pt]
 			\draw[blue,thick](0,-2) node (a)  {$a$};
 			\draw[blue,thick](-1,-3) node (b)  {$b$};
 			\draw[blue,thick](1,-3) node (c)  {$c$};
 			\draw [blue,thick](0, -5) node (u)  {$u$};
 			\draw [blue,thick](-1, -7) node (x)  {$x$};
 			\draw [blue,thick](1, -7) node (y)  {$y$};
 		}
 		\path[draw,thick,black!60!green]
 		(a) edge node[name=la,pos=0.7, above] {\color{blue} $1$} (c)
 		
 		(c) edge node[name=la,pos=0.5, below] {\color{blue}} (u)
 		(b) edge node[name=la,pos=0.5, below] {\color{blue} $\delta$\quad\quad} (u)
 		(u) edge node[name=la,pos=0.6, above] {\color{blue}$\gamma$\quad\quad} (x)
 		(u) edge node[name=la,pos=0.6,above] {\color{blue}  \quad$1$} (y);
 		

 		\draw[dashed, red, line width=0.5mm] (b)--++(140:1cm); 
 		\draw[dashed, red, line width=0.5mm] (x)--++(200:1cm); 
 		\draw[dashed, red, line width=0.5mm] (y)--++(340:1cm);
 		\draw[dashed, red, line width=0.5mm] (a)--++(40:1cm); 
 		\draw[dashed, red, line width=0.5mm] (a)--++(90:1cm); 
 		\draw[dashed, red, line width=0.5mm] (a)--++(140:1cm);

 		\draw[blue] (-1.5, -7.4) node {$\eta$}; 
 		\draw[blue] (1.5, -7.4) node {$\eta$}; 
 		\draw[blue] (-1.2, -2.5) node {$1$};
 		\draw[blue] (0.6, -1.8) node {$\delta$};
 		\draw[blue] (-0.6, -1.8) node {$\gamma$};
 		\draw[blue] (-0.15, -1.5) node {$\eta$};
 		
%
%
%
 		\end{tikzpicture}
 		-	  	\end{center}
 	\caption{Colors on the edges connecting $x$ and $y$ to $b$}
 	\label{f1}
 \end{figure}
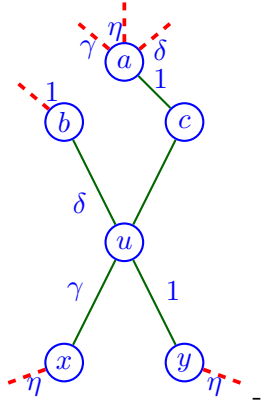

 Let $\varphi(bu)=\delta$.  Again, $\delta\in \pbar(a)$. 
 Figure~\ref{f1} depicts the colors and missing colors  on these specified edges and vertices, respectively.  Clearly, $\delta\ne 1, \gamma$.
 Since $a$ and $b$ are $(1,\delta)$-linked with respect to $\varphi$, $\eta\ne \delta$. Otherwise, $b$ and $y$ would be $(1,\delta)$-linked. 
 Thus, $\gamma, \delta$ and $\eta$ are pairwise distinct.   
 
 The claim below is simple though it is actually the breakthrough for 
 showing Lemma~\ref{lemma:class2-with-fullDpair2}. Without using the symmetry
 of the 4-cycle $abuca$, 
 it seems extremely difficult to prove the same result. 
 
 \begin{CLA}\label{cla:claim2}
 	It holds that $ub\in P_y(\eta,\delta)$ and $P_y(\eta,\delta)$ meets $u$ before $b$. 
 \end{CLA}
 
 \pf Let $\varphi'$ be obtained from $\varphi$ by coloring $ab$ by $\delta$ and uncoloring $bu$. Note that $\pbar'(b)=1, \pbar'(u)=\delta$ and $\varphi'(uy)=1$.
 Thus $F^*=(u, ub,b, uy, y)$ is a multifan and so $u$ and $y$ are $(\eta, \delta)$-linked by Lemma~\ref{thm:vizing-fan1}~\eqref{thm:vizing-fan1b}. By uncoloring $ab$ and coloring 
 $bu$ by $\delta$, we get back the original coloring $\varphi$. Therefore, under the coloring $\varphi$, $u\in P_y(\eta, \delta)$ and $P_y(\eta,\delta)$ meets $u$ before $b$. 
 \qed

 We apply the following operations based on $\varphi$:
 \[
 \begin{bmatrix}
 ux& P_{[u,y]}(\eta,\delta)  &  ub &  P_u(1,\gamma)  & ab\\
 \gamma\rightarrow \eta& \delta/\eta &  \delta \rightarrow 1&
 1/\gamma&  \delta\end{bmatrix}.
 \]
 By Claim~\ref{cla:claim1},  $P_u(1,\gamma)$ does not end at  $b$. 
 In any case, the above operations give  
 an edge $\Delta$-coloring of $G$. This contradicts  the earlier assumption that $\chi'(G)=\Delta+1$.
 
 \medskip 
 
 {\bf \noindent Case 2: $\pbar(x)\ne \pbar(y)$.}
 
 Let 
 $$ \varphi(bu)=\alpha, \quad  \varphi(ux)=\beta,  \quad \pbar(x)=\tau, \quad \text{and} \quad \pbar(y)=\gamma.$$	
 As $\varphi(uy)=\pbar(b)=1$, $1\not\in \{\alpha,\beta,\gamma\}$.
 Also, since $a$ and $b$ are $(1,\alpha)$-linked, $\gamma\ne \alpha$. 
 Otherwise, $b$ and $y$ would be $(1,\alpha)$-linked. 
 Since both $K$ and $K^*$ are Kierstead paths and 
 $\pbar(x)\cup \pbar(y)\subseteq \pbar(a)\cup \pbar(b)$, we have $\alpha,\beta,\tau,\gamma\in\pbar(a)$. 
 \begin{CLA}
 	We may assume $\pbar(x)=\tau=1$. 
 \end{CLA} 
 \pf  If $uy\not\in P_x(1,\tau)$, we simply do a $(1,\tau)$-swap at $x$. 
 Thus, we assume that $u\in P_x(1,\tau)$. We first do a $(1,\tau)$-swap at $b$, then an $(\alpha,\tau)$-swap at $x$. Then we do a $(\gamma,\tau)$-swap at $b$. Finally, a $(1,\gamma)$-swap at $b$ and a $(1,\alpha)$-swap at $x$
 give the desired coloring. 
 \qed 
 
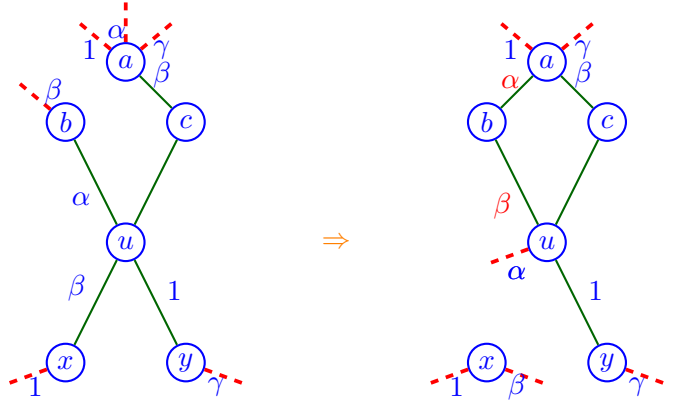
\begin{figure}[!htb]
 	\begin{center}
 		\begin{tikzpicture}[scale=0.8]
 		
 		{\tikzstyle{every node}=[draw ,circle,fill=white, minimum size=0.5cm,
 			inner sep=0pt]
 			\draw[blue,thick](0,-2) node (a)  {$a$};
 			\draw[blue,thick](-1,-3) node (b)  {$b$};
 			\draw[blue,thick](1,-3) node (c)  {$c$};
 			\draw [blue,thick](0, -5) node (u)  {$u$};
 			\draw [blue,thick](-1, -7) node (x)  {$x$};
 			\draw [blue,thick](1, -7) node (y)  {$y$};
 		}
 		\path[draw,thick,black!60!green]
 		(a) edge node[name=la,pos=0.7, above] {\color{blue} $\beta$} (c)
 		
 		(c) edge node[name=la,pos=0.5, below] {\color{blue}} (u)
 		(b) edge node[name=la,pos=0.5, below] {\color{blue} $\alpha$\quad\quad} (u)
 		(u) edge node[name=la,pos=0.6, above] {\color{blue}$\beta$\quad\quad} (x)
 		(u) edge node[name=la,pos=0.6,above] {\color{blue}  \quad$1$} (y);
 		

 		\draw[dashed, red, line width=0.5mm] (b)--++(140:1cm); 
 		\draw[dashed, red, line width=0.5mm] (x)--++(200:1cm); 
 		\draw[dashed, red, line width=0.5mm] (y)--++(340:1cm);
 		\draw[dashed, red, line width=0.5mm] (a)--++(40:1cm); 
 			\draw[dashed, red, line width=0.5mm] (a)--++(90:1cm); 
 		\draw[dashed, red, line width=0.5mm] (a)--++(140:1cm);

 		\draw[blue] (-1.5, -7.4) node {$1$}; 
 		\draw[blue] (1.5, -7.4) node {$\gamma$}; 
 		\draw[blue] (-1.2, -2.5) node {$\beta$};
 		\draw[blue] (0.6, -1.8) node {$\gamma$};
 		\draw[blue] (-0.6, -1.8) node {$1$};
 			\draw[blue] (-0.15, -1.5) node {$\alpha$};

%

 		\draw [orange,thick](3.5, -5) node (t)  {$\Rightarrow$};

 		\begin{scope}[shift={(7,0)}]
 		{\tikzstyle{every node}=[draw ,circle,fill=white, minimum size=0.5cm,
 			inner sep=0pt]
 			\draw[blue,thick](0,-2) node (a)  {$a$};
 			\draw[blue,thick](-1,-3) node (b)  {$b$};
 			\draw[blue,thick](1,-3) node (c)  {$c$};
 			\draw [blue,thick](0, -5) node (u)  {$u$};
 			\draw [blue,thick](-1, -7) node (x)  {$x$};
 			\draw [blue,thick](1, -7) node (y)  {$y$};
 		}
 		\path[draw,thick,black!60!green]
 		(a) edge node[name=la,pos=0.7, above] {\color{blue} $\beta$} (c)
 		(a) edge node[name=la,pos=0.7, above] {\color{red} $\alpha$} (b)
 		
 		(c) edge node[name=la,pos=0.5, below] {\color{blue}} (u)
 		(b) edge node[name=la,pos=0.5, below] {\color{red} $\beta$\quad\quad} (u)
 		(u) edge node[name=la,pos=0.6,above] {\color{blue}  \quad$1$} (y);
 		

 		\draw[dashed, red, line width=0.5mm] (u)--++(200:1cm); 
 		\draw[dashed, red, line width=0.5mm] (x)--++(200:1cm); 
 		\draw[dashed, red, line width=0.5mm] (x)--++(340:1cm);
 		\draw[dashed, red, line width=0.5mm] (y)--++(340:1cm);
 		\draw[dashed, red, line width=0.5mm] (a)--++(40:1cm); 
 		\draw[dashed, red, line width=0.5mm] (a)--++(140:1cm);

 		\draw[blue] (-1.5, -7.4) node {$1$}; 
 		\draw[blue] (-0.5, -7.4) node {$\beta$}; 
 		\draw[blue] (1.5, -7.4) node {$\gamma$}; 
 		\draw[blue] (0.6, -1.8) node {$\gamma$};
 		\draw[blue] (-0.6, -1.8) node {$1$};
 		\draw[blue] (-0.5, -5.5) node {$\alpha$}; 
 			\draw[blue] (-0.5, -5.5) node {$\alpha$};

%
 		\end{scope}
 		\end{tikzpicture}
 		-	  	\end{center}
 	\caption{Colors on the edges connecting $x$ and $y$ to $b$}
 	\label{f2}
 \end{figure}

 Since $ux\in P_x(1,\beta)$, and $a$ and $b$ are $(1,\beta)$-linked, we do a $(1,\beta)$-swap at $b$. 
 Now we color $ab$ by $\alpha$, recolor $bu$ by $\beta$
 and uncolor $ux$, see Figure~\ref{f2} for a depiction. 
 
 Note that 
 $$
 F^*=(u, ux,x,uy,y), \quad K^*=(x,xu,u,ub,b,ba,a)
 $$
 are, respectively,  a multifan and 
 a Kierstead path. 
 By Lemma~\ref{thm:vizing-fan1}~\eqref{thm:vizing-fan1b}, $u$ and $y$ are $(\alpha,\gamma)$-linked, and $u$
 and $x$ are $(\alpha,\beta)$-linked and $(1,\alpha)$-linked. 
 Thus, we do an $(\alpha,\gamma)$-swap at $a$, an $(\alpha,\beta)$-swap at $a$, a $(1,\alpha)$-swap at $a$, 
 and then an $(\alpha,\gamma)$-swap at $a$. Now $P_u(\alpha,\beta)=uba$, contradicting 
 Lemma~\ref{thm:vizing-fan1}~\eqref{thm:vizing-fan1b} that $u$ and $x$ are $(\alpha,\beta)$-linked. 
 The proof  is now completed. 
 \qed


\begin{thebibliography}{10}
 	
 	\bibitem{MR4183149}
 	S. Bonvicini and A. Vietri.
 	\newblock A {M}\"{o}bius-type gluing technique for obtaining edge-critical
 	graphs.
 	\newblock {\em Ars Math. Contemp.}, 19(2):209--229, 2020.
 	
 	\bibitem{MR2028248}
 	S. Gr\"{u}newald and E. Steffen.
 	\newblock Independent sets and 2-factors in edge-chromatic-critical graphs.
 	\newblock {\em J. Graph Theory}, 45(2):113--118, 2004.
 	
 	\bibitem{MR1460574}
 	A.~J.~W. Hilton and C.~Zhao.
 	\newblock Vertex-splitting and chromatic index critical graphs.
 	\newblock volume~76, pages 205--211. 1997.
 	\newblock Second International Colloquium on Graphs and Optimization
 	(Leukerbad, 1994).
 	
 	\bibitem{Holyer}
 	I. Holyer.
 	\newblock The {NP}-completeness of edge-coloring.
 	\newblock {\em SIAM J. Comput.}, 10(4):718--720, 1981.
 	
 	\bibitem{MR1874750}
 	Z. Song.
 	\newblock A further extension of {Y}ap's construction for {$\Delta$}-critical
 	graphs.
 	\newblock {\em Discrete Math.}, 243(1-3):283--290, 2002.
 	
 	\bibitem{StiebSTF-Book}
 	M. Stiebitz, D. Scheide, B. Toft, and L.~M. Favrholdt.
 	\newblock {\em Graph edge coloring}.
 	\newblock Wiley Series in Discrete Mathematics and Optimization. John Wiley \&
 	Sons, Inc., Hoboken, NJ, 2012.
 	\newblock Vizing's theorem and Goldberg's conjecture, With a preface by
 	Stiebitz and Toft.
 	
 	\bibitem{MR3417245}
 	A. Vietri.
 	\newblock An analogy between edge colourings and differentiable manifolds, with
 	a new perspective on 3-critical graphs.
 	\newblock {\em Graphs Combin.}, 31(6):2425--2435, 2015.
 	
 	\bibitem{Vizing64}
 	V.~G. Vizing.
 	\newblock On an estimate of the chromatic class of a {$p$}-graph.
 	\newblock {\em Diskret. Analiz}, (3):25--30, 1964.
 	
 	\bibitem{vizing-2factor}
 	V.~G. Vizing.
 	\newblock The chromatic class of a multigraph.
 	\newblock {\em Kibernetika (Kiev)}, 1965(3):29--39, 1965.
 	
 	\bibitem{Vizing-2-classes}
 	V.~G. Vizing.
 	\newblock Critical graphs with given chromatic class.
 	\newblock {\em Diskret. Analiz No.}, 5:9--17, 1965.
 	
 \end{thebibliography}
\end{document}